\documentclass{article}
\usepackage{amssymb}
\usepackage{amsfonts}
\usepackage{amsmath}

\newtheorem{theorem}{Theorem}[section]

\newtheorem{corollary}[theorem]{Corollary}

\newtheorem{definition}[theorem]{Definition}
\newtheorem{example}[theorem]{Example}

\newtheorem{lemma}[theorem]{Lemma}

\newtheorem{proposition}[theorem]{Proposition}
\newtheorem{remark}[theorem]{Remark}

\numberwithin{equation}{section}
\newenvironment{proof}[1][Proof]{\noindent\textbf{#1.} }{\ \rule{0.5em}{0.5em}}
\input{tcilatex}

\begin{document}

\title{A family of measures associated with iterated function systems}
\author{Palle E. T. Jorgensen\thanks{%
Work supported in part by the NSF.} \\
Department of Mathematics\\
University of Iowa\\
Iowa City, Iowa 52242}
\maketitle

\begin{abstract}
Let $(X,d)$ be a compact metric space, and let an iterated function system
(IFS) be given on $X$, i.e., a finite set of continuous maps $\sigma _{i}$:~$%
X\rightarrow X$, $i=0,1,\cdots ,N-1$. The maps $\sigma _{i}$ transform the
measures $\mu $ on $X$ into new measures $\mu _{i}$. If the diameter of $%
\sigma _{i_{1}}\circ \cdots \circ \sigma _{i_{k}}(X)$ tends to zero as $%
k\rightarrow \infty $, and if $p_{i}>0$ satisfies $\sum_{i}p_{i}=1$, then it
is known that there is a unique Borel probability measure $\mu $ on $X$ such
that 
\begin{equation}
\mu =\dsum_{i}p_{i}~\mu _{i}  \tag{*}
\end{equation}

In this paper, we consider the case when the $p_{i}$s are replaced with a
certain system of sequilinear functionals. This allows us to study the
variable coefficient case of (*), and moreover to understand the analog of
(*) which is needed in the theory of wavelets.
\end{abstract}

\section{Introduction\label{Intro}}

A finite system of continuous functions $\sigma _{i}$:~$X\rightarrow X$ in a
compact metric space $X$ is said to be an \textit{iterated function system}
(IFS) if there is a mapping $\sigma $:~$X\rightarrow X$, onto $X$, such that 
\begin{equation}
\sigma \circ \sigma _{i}=id_{X}  \label{IntroEq1}
\end{equation}

If there is a constant $0<c<1$ such that 
\begin{equation}
d(\sigma _{i}(x),\sigma _{i}(y))\leq c\;d(x,y),\;x,y\in X\text{,}
\label{IntroEq2}
\end{equation}%
then we say that the IFS is \textit{contractive}. In that case, there is,
for every configuration $p_{i}>0$, $\sum_{i}p_{i}=1$, a unique Borel
probability measure $\mu $, $\mu =\mu _{(p)}$ on $X$ such that 
\begin{equation}
\mu =\sum p_{i}\mathcal{\;}\mu \circ \sigma _{i}^{-1}\text{.}
\label{IntroEq3}
\end{equation}

This follows from a theorem of Hutchinson \cite{Hut81}. The mappings $\sigma
_{i}$ might be defined initially on some Euclidean space $E$. If the
contractivity (\ref{IntroEq2}) is assumed, then there is a unique compact
subset $X\subset E$ such that 
\begin{equation}
X=\bigcup_{i}\sigma _{i}(X)\text{,}  \label{IntroEq4}
\end{equation}%
and this set $X$ is the support of $\mu $.

\begin{example}
\label{IntroEx1}Let $E=\mathbb{R}$, $\sigma _{0}(x)=\frac{x}{3}$, $\sigma
_{1}=\frac{x+2}{3}$, $p_{0}=p_{1}=\frac{1}{2}$. In this case, $X$ is the
familiar middle-third Cantor set, and $\mu $ is the Cantor measure supported
in $X$ with Hausdorff dimension $d=\frac{\ln 2}{\ln 3}$. But at the same
time $X$ may be identified with the compact Cartesian product $X\cong \prod 
\mathbb{Z}_{2}=\mathbb{Z}_{2}^{\mathbb{N}}$, where $\mathbb{Z}_{2}=\mathbb{Z}%
\diagup 2\mathbb{Z=}\{0,1\}$, and $\mathbb{N=}\{0,1,2,\cdots \}$, and $\mu $
is the infinite product measure on $\mathbb{Z}_{2}\times \mathbb{Z}%
_{2}\times \cdots $ with weights $(\frac{1}{2},\frac{1}{2})$ on each factor.
\end{example}

\begin{example}
\label{IntroEx2}Let $E=\mathbb{R}$, $\sigma _{0}(x)=\frac{x}{2}$, $\sigma
_{1}(x)=\frac{x+1}{2}$, $p_{0}=p_{1}=\frac{1}{2}$. In this case, $X=[0,1]$,
i.e., the compact unit interval, and $\mu $ is the restriction to $[0,1]$ of
the usual Lebesgue measure $dt$ on $\mathbb{R}$.
\end{example}

Let $\mathbb{T}=\mathbb{R\diagup }2\pi \mathbb{Z}=\{z\in \mathbb{C}\mid
|z|=1\}$ be the usual torus. Let $N\in \mathbb{N},$ $N\geq 2$, and let $%
m_{i} $:~$\mathbb{T}\rightarrow \mathbb{C}$, $i=0,1,\cdots ,N-1$ be a system
of $L^{\infty }$-functions such that the $N\times N$ matrix 
\begin{equation}
\frac{1}{\sqrt{N}}\left( m_{j}\left( ze^{i2\pi \frac{k}{N}}\right) \right)
_{j,k=0}^{N-1}\qquad ,z\in \mathbb{T}  \label{IntroEq5}
\end{equation}%
is unitary. Set 
\begin{equation}
S_{j}f(z)=m_{j}(z)f(z^{N})\qquad ,z\in \mathbb{T}\text{, }f\in L^{2}(\mathbb{%
T})\text{.}  \label{IntroEq6}
\end{equation}

Then it is well known \cite{Jor03b} that the operators $S_{j}$ satisfy the
following two relations, 
\begin{equation}
S_{j}^{\ast }S_{k}=\delta _{j,k}I  \label{IntroEq7}
\end{equation}

\begin{equation}
\sum_{j}S_{j}S_{j}^{\ast }=I  \label{IntroEq8}
\end{equation}%
where $I$ denotes the identity operator in the Hilbert space $\mathcal{H}$:~$%
=L^{2}(\mathbb{T})$. The converse implication also holds, see \cite{BrJo02}.
Systems of isometries satisfying (\ref{IntroEq7}) -- (\ref{IntroEq8}) are
called \textit{representations of the Cuntz algebra} $\mathcal{O}_{N}$, see 
\cite{Cu77}, and the particular representations (\ref{IntroEq6}) are well
known to correspond to multiresolution wavelets; the functions $m_{j}$ are
denoted wavelet filters. These same functions are used in subband filters in
signal processing, see \cite{BrJo02}.

It is easy to see that there is a unique Borel measure $P$ on $[0,1)$ taking
values in the orthogonal projections of $\mathcal{H}$ such that 
\begin{equation}
P\left( \left[ \frac{i_{1}}{N}+\cdots +\frac{i_{k}}{N^{k}},\frac{i_{1}}{N}%
+\cdots +\frac{i_{k}}{N^{k}}+\frac{1}{N^{k}}\right) \right) =S_{i_{1}}\cdots
S_{i_{k}}S_{i_{k}}^{\ast }\cdots S_{i_{1}}^{\ast }\text{.}  \label{IntroEq9}
\end{equation}

Let 
\begin{equation}
e_{n}(z)=z^{n},\;z\in \mathbb{T},n\in \mathbb{Z}\text{.}  \label{IntroEq10}
\end{equation}

\begin{example}
\label{IntroEx3}Let $N\in \mathbb{N,}$ $N\geq 2$, and set $m_{j}(z)=z^{j},$ $%
0\leq j<N$. Then \emph{(}\ref{IntroEq5}\emph{)} is satisfied, and we have 
\begin{equation*}
\left\{ 
\begin{array}{l}
S_{0}^{\ast }e_{0}=e_{0} \\ 
S_{j}^{\ast }e_{0}=0\qquad ,~0<j<N\text{.}%
\end{array}%
\right.
\end{equation*}
\end{example}

It follows easily that the corresponding measure 
\begin{equation}
\mu _{0}(\cdot )\text{:~}=\left\Vert P(\cdot )e_{0}\right\Vert ^{2}
\label{IntroEq11}
\end{equation}%
on $[0,1)$ is the Dirac measure $\delta _{0}$ at $x=0$, i.e., 
\begin{equation}
\delta _{0}(E)=\left\{ 
\begin{array}{c}
1\qquad \text{if }0\in E \\ 
0\qquad \text{if }0\notin E%
\end{array}%
\text{.}\right.  \label{IntroEq12}
\end{equation}

Here $P(\cdot )$ refers to the projection valued measure determined by (\ref%
{IntroEq9}) when the representation of $\mathcal{O}_{N}$ is specified by the
system $m_{j}$:~$=e_{j},$ $0\leq j<N$.

\begin{example}
\label{IntroEx4}Let $N\in \mathbb{N}$, $N\geq 2$, and set 
\begin{equation}
m_{j}(z)=\frac{1}{\sqrt{N}}\sum_{k=0}^{N-1}e^{i2\pi \frac{jk}{N}}~z^{k}\text{%
.}  \label{IntroEq13}
\end{equation}%
Again the condition \emph{(}\ref{IntroEq5}\emph{)} is satisfied, and one
checks that 
\begin{equation}
S_{j}^{\ast }e_{0}=\frac{1}{\sqrt{N}}e_{0}\qquad ,0\leq j<N\text{;}
\label{IntroEq14}
\end{equation}%
and now the measure $\mu _{0}(\cdot )=\left\Vert P(\cdot )e_{0}\right\Vert
^{2}$ on $[0,1)$ is the restriction to $[0,1)$ of the Lebesgue measure on $%
\mathbb{R}$. It is well known that the wavelet corresponding to \emph{(}\ref%
{IntroEq13}\emph{)} is the familiar Haar wavelet corresponding to $N$-adic
subdivision, see \cite{BrJo02}. It is also known that generally, for
wavelets other than the Haar systems, the corresponding representations 
\emph{(}\ref{IntroEq6}\emph{)} of $\mathcal{O}_{N}$ does not admit a
simultaneous eigenvector $f$, i.e., there is no solution $f\in \mathcal{H}%
\diagdown \{0\},$ $\lambda _{j}\in \mathbb{C}$ to the joint eigenvalue
problem%
\begin{equation}
S_{j}^{\ast }f=\lambda _{j}f\qquad ,0\leq j<N\text{.}  \label{IntroEq15}
\end{equation}
\end{example}

\begin{proposition}
\label{IntroProp1}Let $N\in \mathbb{N},$ $N\geq 2$, and let $(S_{j})_{0\leq
j<N}$ be a representation of $\mathcal{O}_{N}$ on a Hilbert space $\mathcal{H%
}$. Suppose there is a solution $f\in \mathcal{H},$ $\left\Vert f\right\Vert
=1$, to the eigenvalue problem \emph{(}\ref{IntroEq15}\emph{)} for some $%
\lambda _{j}\in \mathbb{C}$. Then $\sum \left\vert \lambda _{j}\right\vert
^{2}=1$, and the measure $\mu $\emph{:}~$=\left\Vert P(\cdot )f\right\Vert
^{2}$ satisfies 
\begin{equation}
\mu =\sum_{j=0}^{N-1}\left\vert \lambda _{j}\right\vert ^{2}\mu \circ \sigma
_{j}^{-1}  \label{IntroEq16}
\end{equation}%
where $\sigma _{j}(x)=\frac{x+j}{N},$ $\mu \circ \sigma _{j}^{-1}(E)=\mu
(\sigma _{j}^{-1}(E))$ for Borel sets $E\subset \lbrack 0,1)$, and $\sigma
_{j}^{-1}(E)$\emph{:}~$=\{x\mid \sigma _{j}(x)\in E\}$.
\end{proposition}

\begin{proof}
The reader may prove the proposition directly from the definitions, but the
conclusion may also be obtained as a special case of the theorem in the next
section.
\end{proof}

\section{Measures And Iterated Function Systems\label{Measures}}

Let $(X,d)$ be a compact metric space, and let $(\sigma _{j})_{0\leq j<N}$
be an $N$-adic iterated function system (IFS). We say that the system is 
\textit{complete} if 
\begin{equation}
\lim_{k\rightarrow \infty }\text{diameter}~\left( \sigma _{i_{1}}\circ
\cdots \circ \sigma _{i_{k}}(X)\right) =0\text{, }  \label{MeasEq1}
\end{equation}

We say that the IFS is \textit{non-overlapping} if for each $k$ the sets 
\begin{equation}
A_{k}(i_{1},\cdots ,i_{k})\text{:~}=\sigma _{i_{1}}\circ \cdots \circ \sigma
_{i_{k}}(X)  \label{MeasEq2}
\end{equation}%
are disjoint, i.e., for every $k$, the sets $A_{k}(i_{1},\cdots ,i_{k})$ are
mutually disjoint for different multi-indices, i.e., different points in 
\begin{equation*}
\underset{k\text{-times}}{\underbrace{\mathbb{Z}_{N}\times \cdots \times 
\mathbb{Z}_{N}}}
\end{equation*}%
where $\mathbb{Z}_{N}$:~$=\{0,1,\cdots ,N-1\}$.

\begin{remark}
\label{MeasRem}It is immediate that, if a given IFS $(\sigma _{j})_{0\leq
j<N}$ arises as a system of \textit{distinct} branches of the inverse of a
single mapping $\sigma $\emph{:}~$X\rightarrow X$, i.e., if $\sigma (\sigma
_{i}(x))=x$ for $x\in X$, and $0\leq i<N$, then the partition system $\sigma
_{i_{1}}\circ \cdots \circ \sigma _{i_{k}}(X)$ is non-overlapping.
\end{remark}

\begin{theorem}
\label{MeasTheo1}Let $N\in \mathbb{N},$ $N\geq 2$, and a Hilbert space $%
\mathcal{H}$. Let $(\sigma _{j})_{0\leq j<N}$ be an IFS which is complete
and non-overlapping. Then there is a unique projection valued measure $P$
defined on the Borel subsets of $X$ such that 
\begin{equation}
P\left( A_{k}(i_{1},\cdots ,i_{k})\right) =S_{i_{1}}\cdots
S_{i_{k}}S_{i_{k}}^{\ast }\cdots S_{i_{1}}^{\ast }\text{.}  \label{MeasEq3}
\end{equation}%
This measure satisfies$\emph{:}$

\noindent \emph{(}a\emph{)} $P(E)=P(E)^{\ast }=P(E)^{2}$, $E\in \mathcal{B}%
(X)=$ the Borel subsets of $X$.

\noindent \emph{(}b\emph{)} $\int_{X}dP(x)=I_{\mathcal{H}}$

\noindent \emph{(}c\emph{)} $P(E)P(F)=0$ if $E,F\in \mathcal{B}(X)$ and $%
E\cap F=\varnothing $.

\noindent \emph{(}d\emph{)} $\sum_{j=0}^{N-1}S_{j}P(\sigma
_{j}^{-1}(E))S_{j}^{\ast }=P(E)$, $E\in \mathcal{B}(X)$.

It follows in particular that, for every $f\in \mathcal{H}$, the measure $%
\mu _{f}(\cdot )$:$=\left\Vert P(\cdot )f\right\Vert ^{2}$ satisfies 
\begin{equation}
\sum_{j=0}^{N-1}\mu _{S_{j}^{\ast }f}\circ \sigma _{j}^{-1}=\mu _{f}\text{,}
\label{MeasEq4}
\end{equation}%
or equivalently 
\begin{equation}
\sum_{j=0}^{N-1}\int_{X}\psi \circ \sigma _{j}\mathcal{\;}d\mu _{S_{j}^{\ast
}f}=\int_{X}\psi \mathcal{\;}d\mu _{f}  \label{MeasEq5}
\end{equation}%
for all bounded Borel functions $\psi $ on $X$.
\end{theorem}

\begin{corollary}
\label{MeasCor1}Let N$\in \mathbb{N}$, $N\geq 2$, be given, and consider a
representation $(S_{j})_{0\leq j<N}$ of $\mathcal{O}_{N}$, and an associated
IFS which is complete and non-overlapping. Let $P(\cdot )$ be the
corresponding projection valued measure, i.e., 
\begin{equation}
P\left( \sigma _{i_{1}}\circ \cdots \circ \sigma _{i_{k}}(X)\right)
=S_{i_{1}}\cdots S_{i_{k}}S_{i_{k}}^{\ast }\cdots S_{i_{1}}^{\ast }\text{.}
\label{MeasEq6}
\end{equation}%
For $f\in \mathcal{H},$ $\left\Vert f\right\Vert =1$, set $\mu _{f}(\cdot )$%
\emph{:}$=\left\Vert P(\cdot )f\right\Vert ^{2}$. Let $\mathfrak{A}$ be the
abelian $C^{\ast }$-algebra generated by the projections $S_{i_{1}}\cdots
S_{i_{k}}S_{i_{k}}^{\ast }\cdots S_{i_{1}}^{\ast }$ and let $\mathcal{H}_{f}$
be the closure of $\mathfrak{A}f$. Then there is a unique isometry $V_{f}$:~$%
L^{2}(\mu _{f})\rightarrow \mathcal{H}_{f}$ of $L^{2}(\mu _{f})$ onto $%
\mathcal{H}_{f}$ such that 
\begin{equation}
V_{f}(1)=f\text{,}  \label{MeasEq7}
\end{equation}%
and 
\begin{equation}
V_{f}M_{\chi _{A_{k}(i)}}V_{f}^{\ast }=S_{i_{1}}\cdots
S_{i_{k}}S_{i_{k}}^{\ast }\cdots S_{i_{1}}^{\ast }  \label{MeasEq8}
\end{equation}%
where $M_{\chi _{A_{k}(i)}}$ is the operator which multiplies by the
indicator function of 
\begin{equation}
A_{k}(i)\emph{:}=\sigma _{i_{1}}\circ \cdots \circ \sigma _{i_{k}}(X)\text{.}
\label{MeasEq9}
\end{equation}
\end{corollary}

\begin{proof}
(Theorem \ref{MeasTheo1}) We refer to the paper \cite{Jor03} for a more
complete discussion. With the assumptions, we note that for every $k$, and
every multi-index $i=(i_{1},\cdots ,i_{k})$ we have an abelian algebra of
functions $\mathcal{F}_{k}$ spanned by the indicator functions $\chi
_{A_{k}(i)}$ where $A_{k}(i)$:$=\sigma _{i_{1}}\circ \cdots \circ \sigma
_{i_{k}}(X)$. Since, for every $k$, we have the non-overlapping unions 
\begin{equation}
\bigcup_{i}A_{k+1}(i_{1},i_{2},\cdots ,i_{k},i)=A_{k}(i_{1},\cdots ,i_{k})%
\text{,}  \label{MeasEq10}
\end{equation}%
there is a natural embedding $\mathcal{F}_{k}\subset \mathcal{F}_{k+1\text{.}%
}$ We wish to define the projection valued measure $P$ as an operator valued
map on functions on $X$ in such a way that $\int_{X}\chi
_{_{A_{k}(i)}}dP=S_{i_{1}}\cdots S_{i_{k}}S_{i_{k}}^{\ast }\cdots
S_{i_{1}}^{\ast }$. This is possible since the projections $S_{i_{1}}\cdots
S_{i_{k}}S_{i_{k}}^{\ast }\cdots S_{i_{1}}^{\ast }$ are mutually orthogonal
when $k$ is fixed, and $(i_{1},\cdots ,i_{k})$ varies over $(\mathbb{Z}%
_{N})^{k}$. In view of the inclusions 
\begin{equation}
\mathcal{F}_{k}\subset \mathcal{F}_{k+1}\text{,}  \label{MeasEq11}
\end{equation}%
it follows that $\bigcup_{k}\mathcal{F}_{k}$ is an algebra of functions on $%
X $. Since the $N$-adic subdivision system $\{A_{k}(i)\mid k\in \mathbb{N}$, 
$i\in (\mathbb{Z}_{N})^{k}\}$ is complete, it follows that every continuous
function on $X$ is the uniform limit of a sequence of functions in $%
\bigcup_{k}\mathcal{F}_{k}$. Using now a standard extension procedure for
measures, we conclude that the projection valued measure $P(\cdot )$ exists,
and that it has the properties listed in the theorem. The reader is referred
to \cite{Nel70} for additional details on the extension from $\bigcup_{k}%
\mathcal{F}_{k}$ to the Borel function on $X$.
\end{proof}

\begin{proof}[Proof of Corollary \protect\ref{MeasCor1}]
Let the systems $(S_{j})$ and $(\sigma _{j})$ be as in the statement of the
theorem. To define $V_{f}$:~$L^{2}(\mu _{f})\rightarrow \mathcal{H}_{f}$, we
set 
\begin{equation}
V_{f}(1)=f,\text{and }V_{f}\chi _{A_{k}(i)}=S_{i_{1}}\cdots
S_{i_{k}}S_{i_{k}}^{\ast }\cdots S_{i_{1}}^{\ast }f\text{.}  \label{MeasEq12}
\end{equation}%
It is clear from the theorem that $V_{f}$ defined this way extends to an
isometry of $L^{2}(\mu _{f})$ onto $\mathcal{H}_{f}$, and a direct
verification reveals that the covariance relation (\ref{MeasEq8}) is
satisfied.

It remains to prove (d) in theorem \ref{MeasTheo1}, or equivalently to prove
(\ref{MeasEq4}) for every $f\in \mathcal{H}$. The argument is based on the
same approximation procedure as we used above, starting with the algebra $%
\bigcup_{k}\mathcal{F}_{k}$. Note that 
\begin{eqnarray*}
&&\int_{X}\chi _{_{A_{k}(\alpha _{1},\cdots ,\alpha _{k})}}(\sigma _{i}(x))%
\mathcal{\;}d\mu _{S_{i}^{\ast }f}(x) \\
&=&\delta _{i,\alpha _{1}}\int_{X}\chi _{_{A_{k}(\alpha _{1},\cdots ,\alpha
_{k})}}(x)\mathcal{\;}d\mu _{S_{i}^{\ast }f}(x) \\
&=&\delta _{i,\alpha _{1}}\left\Vert S_{\alpha _{k}}^{\ast }\cdots S_{\alpha
_{2}}^{\ast }S_{i}^{\ast }f\right\Vert ^{2} \\
&=&\delta _{i,\alpha _{1}}\int_{X}\chi _{A_{k}(\alpha )}(x)\mathcal{\;}d\mu
_{f}(x)\text{.}
\end{eqnarray*}%
Summing over $i$, we get 
\begin{equation*}
\sum_{i}\int_{X}\chi _{A_{k}(\alpha )}\circ \sigma _{i}\;d\mu _{S_{i}^{\ast
}f}=\int_{X}\chi _{A_{k}(\alpha )}d\mu _{f}=\left\Vert S_{\alpha _{k}}^{\ast
}\cdots S_{\alpha _{1}}^{\ast }{}_{i}f\right\Vert ^{2}\text{.}
\end{equation*}%
The desired identity (\ref{MeasEq5}) now follows by yet another application
of the standard approximation argument which was used in the proof of the
first part of the theorem.
\end{proof}

The simplest subdivision system is the one where the subdivisions are given
by the $N$-adic fractions $\frac{\alpha _{1}}{N}+\frac{\alpha _{2}}{N^{2}}%
+\cdots +\frac{\alpha _{k}}{N^{k}}$ where $\alpha _{i}\in \mathbb{Z}%
_{N}=\{0,1,\cdots ,N-1\}$. Setting $\sigma _{j}(x)=\frac{x+j}{N}$, $j\in 
\mathbb{Z}_{N}$, we note that 
\begin{equation*}
\sigma _{\alpha _{1}}\circ \cdots \circ \sigma _{\alpha _{k}}([0,1))=\left[ 
\frac{\alpha _{1}}{N}+\cdots +\frac{\alpha _{k}}{N^{k}},\frac{\alpha _{1}}{N}%
+\cdots +\frac{\alpha _{k}}{N^{k}}+\frac{1}{N^{k}}\right) \text{.}
\end{equation*}%
As a result both the projection valued measure $P(\cdot )$ and the
individual measures $\mu _{f}(\cdot )=\left\Vert P(\cdot )f\right\Vert ^{2}$
are defined on the Borel subsets of $[0,1)$. If $\hat{P}(t)$:$%
=\int_{0}^{1}e^{itx}dP(x)$, then $\left\langle f\mid \hat{P}%
(t)f\right\rangle =\hat{\mu}_{f}(t)$ is the usual Fourier transform of the
measure $\mu _{f}$ for $f\in \mathcal{H}$.

Moreover, by the Spectral theorem \cite{Nel70}, there is a selfadjoint
operator $D$ with spectrum contained in $[0,1]$ such that 
\begin{equation*}
\widehat{P}(t)=e^{itD}\text{,}
\end{equation*}%
see \cite{Nel70}. In fact, the spectrum of $D$ is equal to the support of
the projection valued measure $P(\cdot )$.

\begin{corollary}
\label{MeasCor2}Suppose the $N$-adic partition system used in the theorem is
given by the $N$-adic fractions as in \ref{MeasEq10}$.$ Then the Fourier
transform 
\begin{equation}
\hat{\mu}_{f}(t)=\int_{0}^{1}e^{itx}\;d\mu _{f}(x)\qquad ,t\in \mathbb{R}
\label{MeasEq13}
\end{equation}%
of the measure $\mu _{f}(\cdot )=\left\Vert P(\cdot )f\right\Vert ^{2}$
satisfies 
\begin{equation}
\hat{\mu}_{f}(t)=\sum_{k=0}^{N-1}e^{i\frac{kt}{N}}\mathcal{\;}\widehat{\mu
_{S_{k}^{\ast }f}}(t/N)\text{.}  \label{MeasEq14}
\end{equation}
\end{corollary}

\begin{proof}
With the $N$-adic subdivisions of the unit interval, the maps $\sigma _{k}$
are $\sigma _{k}(x)=\frac{x+k}{N}$ for $k\in \mathbb{Z}_{N}=\{0,1,\cdots
,N-1\}$. Setting $\psi _{t}(x)=e^{itx}$ in (\ref{MeasEq4}), the desired
result (\ref{MeasEq14}) follows immediately.
\end{proof}

In the next result we show that for every $k\in \mathbb{N}$, there is an
approximation formula for the Fourier transform $\hat{\mu}_{f}$ of the
measure $\mu _{f}(\cdot )=\left\Vert P(\cdot )f\right\Vert ^{2}$ involving
the numbers $\left\Vert S_{\alpha _{k}}^{\ast }\cdots S_{\alpha _{1}}^{\ast
}f\right\Vert ^{2}$ as the multi-index $\alpha =(\alpha _{1},\cdots ,\alpha
_{k})$ ranges over $(\mathbb{Z}_{N})^{k}$.

\begin{corollary}
\label{MeasCor3}Let $N\in \mathbb{N}$, $N\geq 2$ be given. Let $%
(S_{j})_{0\leq j<N}$ be a representation of $\mathcal{O}_{N}$ on a Hilbert
space $\mathcal{H}$, and let $P(\cdot )$ be the corresponding
projection-valued measure defined on $\mathcal{B}([0,1))$. Let $f\in 
\mathcal{H},$ $\left\Vert f\right\Vert =1$, and set $\mu _{f}(\cdot
)=\left\Vert P(\cdot )f\right\Vert ^{2}$. Then, for every $k$, we have the
approximation 
\begin{equation}
\left\vert \hat{\mu}_{f}(t)-\sum_{\alpha _{1},\cdots ,\alpha
_{k}}e^{it\left( \frac{\alpha _{1}}{N}+\cdots +\frac{\alpha _{k}}{N^{k}}%
\right) }\left\Vert S_{\alpha }^{\ast }f\right\Vert ^{2}\right\vert \leq
\left\vert t\right\vert N^{-k}  \label{MeasEq15}
\end{equation}%
where the summation is over multi-indices from $(\mathbb{Z}_{N})^{k}$, and $%
S_{\alpha }^{\ast }$\emph{:}$=S_{\alpha _{k}}^{\ast }\cdots S_{\alpha
_{1}}^{\ast }$.
\end{corollary}

\begin{proof}
A $k$-fold iteration of formula (\ref{MeasEq14}) from the previous corollary
yields, 
\begin{equation*}
\hat{\mu}_{f}(t)=\sum_{\alpha _{1},\cdots ,\alpha _{k}}e^{it\left( \frac{%
\alpha _{1}}{N}+\cdots +\frac{\alpha _{k}}{N^{k}}\right) }\widehat{\mu
_{S_{\alpha }^{\ast }f}}(t/N^{k})
\end{equation*}%
and 
\begin{equation*}
\widehat{\mu _{S_{\alpha }^{\ast }f}}(t/N^{k})-\left\Vert S_{\alpha }^{\ast
}f\right\Vert ^{2}=\int_{0}^{1}(e^{itN^{-k}x}-1)\mathcal{\;}d\mu _{S_{\alpha
}^{\ast }f}(x)\text{;}
\end{equation*}%
and therefore 
\begin{eqnarray*}
&&\left\vert \hat{\mu}_{S_{\alpha }^{\ast }f}(tN^{-k})-\left\Vert S_{\alpha
}^{\ast }f\right\Vert ^{2}\right\vert \\
&\leq &\left\vert t\right\vert N^{-k}\int_{0}^{1}x\mathcal{\;}d\mu
_{S_{\alpha }^{\ast }f}(x) \\
&\leq &\left\vert t\right\vert N^{-k}\int_{0}^{1}d\mu _{S_{\alpha }^{\ast
}f}(x) \\
&=&\left\vert t\right\vert N^{-k}\left\Vert S_{\alpha }^{\ast }f\right\Vert
^{2}\text{.}
\end{eqnarray*}%
It follows that the difference on the left-hand side in (\ref{MeasEq15}) is
estimated above in absolute value by 
\begin{eqnarray*}
&&\sum_{\alpha _{1},\cdots ,\alpha _{k}}\left\vert e^{it\left( \frac{\alpha
_{1}}{N}+\cdots +\frac{\alpha _{k}}{N^{k}}\right) }\right\vert \left\vert
t\right\vert N^{-k}\left\Vert S_{\alpha }^{\ast }f\right\Vert ^{2} \\
&=&\left\vert t\right\vert N^{-k}\sum_{\alpha _{1},\cdots ,\alpha
_{k}}\left\Vert S_{\alpha }^{\ast }f\right\Vert ^{2} \\
&=&\left\vert t\right\vert N^{-k}\left\langle f\mid \sum_{\alpha _{1},\cdots
,\alpha _{k}}S_{\alpha }S_{\alpha }^{\ast }f\right\rangle \\
&=&\left\vert t\right\vert N^{-k}\left\langle f\mid f\right\rangle \\
&=&\left\vert t\right\vert N^{-k}\left\Vert f\right\Vert ^{2} \\
&=&\left\vert t\right\vert N^{-k}\text{.}
\end{eqnarray*}
\end{proof}

\begin{definition}
\label{MeasDef1}Let $k\in \mathbb{N}$, and set 
\begin{equation}
x_{k}(\alpha )\emph{:}=\frac{\alpha _{1}}{N}+\cdots +\frac{\alpha _{k}}{N^{k}%
}\;\text{ for }\alpha _{i}\in \{0,1,\cdots ,N-1\}\text{.}  \label{MeasEq16}
\end{equation}%
Let $(S_{i})$ and $(\sigma _{i})$ be as in Corollary \ref{MeasCor2}, and let 
$f\in \mathcal{H},$ $\left\Vert f\right\Vert =1$. We set 
\begin{equation}
\mu _{f}^{(k)}=\sum_{\alpha _{1},\cdots ,\alpha _{k}}\left\Vert S_{\alpha
}^{\ast }f\right\Vert ^{2}\delta _{x_{k}(\alpha )}\text{.}  \label{MeasEq17}
\end{equation}
\end{definition}

These measures form the sequence of measures which we use in the Riemann sum
approximation of Corollary \ref{MeasCor3}; and we are still viewing the
measures $\mu _{f}$ and $\mu _{f}^{(k)}$ as measures on the unit-interval $%
[0,1)$.

\begin{corollary}
\label{MeasCor4} Let $N\in \mathbb{N},$ $N\geq 2$. Let $(S_{i})$ and $%
(\sigma _{j})$ be as in corollary \ref{MeasCor2}, i.e., $(S_{i})$ is in $%
\limfunc{Rep}(\mathcal{O}_{N},\mathcal{H})$ for some Hilbert space $\mathcal{%
H}$, and $\sigma _{j}(x)=\frac{x+j}{N}$ for $x\in \lbrack 0,1)$ and $j\in
\{0,1,\cdots ,N-1\}$. Let $\psi $ be a continuous function on $[0,1)$, and
let $k\in \mathbb{N}$. Then 
\begin{equation}
\left\vert \int_{0}^{1}\psi \mathcal{\;}d\mu _{f}-\int_{0}^{1}\psi \mathcal{%
\;}d\mu _{f}^{(k)}\right\vert \leq N^{-k}\int_{\mathbb{R}}\left\vert t%
\widehat{\psi }(t)\right\vert \mathcal{\;}dt  \label{MeasEq18}
\end{equation}%
where 
\begin{equation}
\widehat{\psi }(t)=\int_{0}^{1}\psi (x)e^{-itx}\mathcal{\;}dx
\label{MeasEq19}
\end{equation}%
is the usual Fourier transform; and we are assuming further that 
\begin{equation*}
\int_{\mathbb{R}}\left\vert t\widehat{\psi }(t)\right\vert \mathcal{\;}%
dt<\infty \text{.}
\end{equation*}
\end{corollary}

\begin{proof}
By the Fourier inversion formula, $\psi (x)=\frac{1}{2\pi }\int_{\mathbb{R}}%
\widehat{\psi }(t)e^{itx}\mathcal{\;}dt$; and we get the following formula
by a change of variables, and by the use of Fubini's theorem: 
\begin{equation}
\int \psi \mathcal{\;}d\mu _{f}-\int \psi \mathcal{\;}d\mu _{f}^{(k)}=\frac{1%
}{2\pi }\int_{\mathbb{R}}\widehat{\psi }(t)\left( \widehat{\mu }_{f}(t)-%
\widehat{\mu _{f}^{(k)}}(t)\right) \mathcal{\;}dt\text{.}  \label{MeasEq20}
\end{equation}%
Since 
\begin{equation}
\widehat{\mu _{f}^{(k)}}(t)=\sum_{\alpha _{1},\cdots ,\alpha
_{k}}e^{itx_{k}(\alpha )}\left\Vert S_{\alpha }^{\ast }f\right\Vert ^{2}%
\text{,}  \label{MeasEq21}
\end{equation}%
the estimate (\ref{MeasEq18}) from Corollary \ref{MeasCor2} applies. An
estimation of the differences in (\ref{MeasEq20}) now yields: 
\begin{equation*}
\left\vert \int \psi \mathcal{\;}d\mu _{f}-\int \psi \mathcal{\;}d\mu
_{f}^{(k)}\right\vert \leq \frac{1}{2\pi }\int_{\mathbb{R}}\left\vert 
\widehat{\psi }(t)\right\vert \left\vert t\right\vert \cdot N^{-k}\mathcal{\;%
}dt
\end{equation*}%
which is the desired result.
\end{proof}

\begin{remark}
\label{MeasRem1}In general, a sequence of probability measures on a compact
Hausdorff space $X$, $(\mu _{k})$ is said to \textit{converge weakly} to the
limit $\mu $ if 
\begin{equation}
\lim_{k\rightarrow \infty }\int_{X}\psi \mathcal{\;}d\mu _{k}=\int_{X}\psi 
\mathcal{\;}d\mu \text{ for all }\psi \in C(X)\text{.}  \label{MeasEq22}
\end{equation}%
However, the conclusion of Corollary \ref{MeasCor5} for the convergence $%
\lim_{k\rightarrow \infty }\mu _{f}^{(k)}=\mu _{f}$ is in fact stronger than
weak convergence as we will show. The notion of weak convergence of measures
is significant in probability theory, see e.g., \cite{Bil71}.
\end{remark}

Since the measures $\mu _{f}$ and $\mu _{f}^{(k)}$ are defined on $\mathcal{B%
}([0,1))$, the corresponding distribution functions $F_{f}$ and $F_{f}^{(k)}$
are defined for $x\in \lbrack 0,1)$ as follows 
\begin{equation}
F_{f}(x)=\mu _{f}([0,x])\text{ and }F_{f}^{(k)}(x)=\mu _{f}^{(k)}([0,x])%
\text{.}  \label{MeasEq23}
\end{equation}

\begin{corollary}
\label{MeasCor5}Let $(S_{i})$ and $(\sigma _{i})$ be as in Corollary \ref%
{MeasCor4}. Let $f\in \mathcal{H}$, $\left\Vert f\right\Vert =1$, be given,
and let $\mu _{f}$ resp., $\mu _{f}^{(k)}$ be the corresponding measures,
with distribution functions $F_{f}$ and $F_{f}^{(k)}$, respectively. Then 
\begin{equation}
\lim_{k\rightarrow \infty }F_{f}^{(k)}(x)=F_{f}(x)\text{.}  \label{MeasEq24}
\end{equation}
\end{corollary}

\begin{proof}
We already proved that the sequence of measures $\mu _{f}^{(k)}$ converges
weakly to $\mu _{f}$ as $k\rightarrow \infty $. Furthermore, it is known
that weak convergence $\mu _{f}^{(k)}\rightarrow \mu _{f}$ implies that (\ref%
{MeasEq24}) holds whenever $x$ is a point of continuity for $F_{f}$, see 
\cite[Theorem 2.3, p.5]{Bil71}. (For the case of the wavelet
representations, it is known that every $x$ is a point of continuity, but
Example \ref{IntroEx3} shows that the measures $\mu _{f}$ are not continuous
in general.) The argument from the proof of Corollary \ref{MeasCor2} shows
that, in general, the points of discontinuity of $F_{f}(\cdot )$ must lie in
the set (\ref{MeasEq16}) of $N$-adic fractions. Using Theorem \ref{MeasTheo1}
and formula (\ref{MeasEq17}), we conclude that if $x_{k}(\alpha )$ is a
point of discontinuity of $F_{f}(\cdot )$, then $\left\vert
F_{f}^{(k+n)}(x_{k}(\alpha ))-F_{f}(x_{k}(\alpha ))\right\vert \leq N^{-k-n}$%
, and therefore 
\begin{equation*}
\lim_{n\rightarrow \infty }F_{f}^{(k+n)}(x_{k}(\alpha ))=F_{f}(x_{k}(\alpha
))
\end{equation*}%
which is the desired conclusion, see (\ref{MeasEq24}).
\end{proof}

\section{The Measures $\protect\mu _{f}$\label{TheMeas}}

In the previous section, we showed that the decomposition theory for
representations of the Cuntz algebra $\mathcal{O}_{N}$ may be analyzed by
the use of projection valued measures on a class of iterated function
systems (IFS). It is known that there is a simple $C^{\ast }$-algebra $%
\mathcal{O}_{N}$ for each $N\in \mathbb{N}$, $N\geq 2$, such that the
representations of $\mathcal{O}_{N}$ are in on-one correspondence with
systems of isometries $(S_{i})$ which satisfy the two relations (\ref%
{IntroEq7})--(\ref{IntroEq8}). The $C^{\ast }$-algebra $\mathcal{O}_{N}$ $is$
defined abstractly on $N$ generators $s_{i}$ which satisfy 
\begin{equation}
\sum_{i=0}^{N-1}s_{i}s_{i}^{\ast }=1\text{ and }s_{i}^{\ast }s_{j}=\delta
_{i,j}1\text{.}  \label{TheMeasEq1}
\end{equation}%
A representation $\rho $ of $\mathcal{O}_{N}$ on a Hilbert space $\mathcal{H}
$ is a $\ast $-homomorphism from $\mathcal{O}_{N}$ into $B(\mathcal{H})=$
the algebra of all bounded operators on $\mathcal{H}$. The set of
representations acting on $\mathcal{H}$ is denoted $\limfunc{Rep}(\mathcal{O}%
_{N},\mathcal{H})$. The connection between $\rho $ and the corresponding $%
(S_{i})$-system is fixed by $\rho (s_{i})=S_{i}$. While the subalgebra $%
\mathcal{C}$ in $\mathcal{O}_{N}$ generated by the monomials $%
s_{i_{j}}\cdots s_{i_{k}}s_{i_{k}}^{\ast }\cdots s_{i_{1}}^{\ast }$ is
maximal abelian in $\mathcal{O}_{N}$, the von Neumann algebra $\mathfrak{A}$
generated by $\rho (\mathcal{C})$ may not be maximally abelian in $B(%
\mathcal{H})$. Whether it is, or not, depends on the representation. It is
known to be maximally abelian if the operators $S_{i}=\rho (s_{i})$ are
given by (\ref{IntroEq6}), and if the functions $m_{i}$ satisfy the usual
subband conditions from wavelet theory. For details, see \cite{Jor01} and 
\cite{BrJo02}. For these representations, $\mathcal{H}=L^{2}(\mathbb{T})$;
and the representations define wavelets 
\begin{equation}
\psi _{i,j,k}(x)\text{:~}=N^{\frac{j}{2}}\psi _{i}(N^{j}x-k)\text{, }%
i=1,\cdots ,N-1\text{, ~}j,k\in \mathbb{Z}  \label{TheMeasEq2}
\end{equation}%
in $L^{2}(\mathbb{R})$. Such wavelets are specified by the functions $\psi
_{1},\cdots ,\psi _{N-1}\in L^{2}(\mathbb{R})$. These representations $\rho $
are called \textit{wavelet representations.} The assertion is that, if $%
\mathfrak{A}$ is defined by a wavelet representation, then $\mathfrak{A}$ is
maximally abelian in $B(L^{2}(\mathbb{T}))$. The operators commuting with $%
\mathfrak{A}$ are denoted $\mathfrak{A}^{\prime }$, and it is easy to see
that $\mathfrak{A}$ is maximally abelian if and only if $\mathfrak{A}%
^{\prime }$ is abelian. An abelian von Neumann algebra $\mathfrak{A}\subset
B(\mathcal{H})$ is said to have a cyclic vector $f$ if the closure of $%
\mathfrak{A}f$ is $\mathcal{H}.$ For $f\in \mathcal{H}$, the closure of $%
\mathfrak{A}f$ is denoted $\mathcal{H}_{f}$. It is known that $\mathfrak{A}$
has a cyclic vector if and only if it is maximally abelian. Clearly, if $f$
is a cyclic vector, then the measure $\mu _{f}(\cdot )$:$~=\left\Vert
P(\cdot )f\right\Vert ^{2}$ determines the other measures $\{\mu _{g}\mid
g\in \mathcal{H}\}$.

\begin{lemma}
\label{TheMeasLem1}Let $\mathfrak{A}\subset B(\mathcal{H})$ be an abelian $%
C^{\ast }$-algebra, and let $\rho $\emph{:~}$C(X)\cong \mathfrak{A}$ be the
Gelfand representation, $X$ a compact Hausdorff space. Let $f\in \mathcal{H}$%
, $\left\Vert f\right\Vert =1$.

\noindent \emph{(}a\emph{)} Then there is a unique Borel measure $\mu $ on $%
X $, and an isometry $V_{f}$\emph{:~}$L^{2}(\mu )\rightarrow \mathcal{H}$,
such that 
\begin{equation}
V_{f}(1)=f\text{,}  \label{TheMeasEq3}
\end{equation}%
\begin{equation}
V_{f}(\psi )=\rho (\psi )f\quad ,\psi \in C(X)\text{,}  \label{TheMeasEq4}
\end{equation}%
and 
\begin{equation}
V_{f}(L^{2}(\mu ))=\mathcal{H}_{f}\text{.}  \label{TheMeasEq5}
\end{equation}

\noindent \emph{(}b\emph{)} Let $f_{i}\in \mathcal{H}$, $\left\Vert
f_{i}\right\Vert =1$, $i=1,2$, and suppose $\mu _{1}<<\mu _{2}$. Setting $k=%
\frac{d\mu _{1}}{d\mu _{2}}$ where $\mu _{i}$\emph{:~}$=\mu _{f_{i}}$, $%
i=1,2 $, then $U\psi =\sqrt{k}\psi $ defines an isometry $U$\emph{:}$%
~L^{2}(\mu _{1})\rightarrow L^{2}(\mu _{2})$, and $W$\emph{:}~$%
=V_{f_{2}}UV_{f_{1}}^{\ast }$\emph{:}~$\mathcal{H}_{f_{1}}\rightarrow 
\mathcal{H}_{f_{2}}$ is in the commutant of $\mathfrak{A}$.
\end{lemma}

\begin{proof}
Part (a) follows from the spectral theorem applied to abelian $C^{\ast }$%
-algebras, see e.g., \cite{Nel70}. To prove (b), let $f_{i}$ be the two
vectors in $\mathcal{H}$, and set $\mu _{i}$:~$=\mu _{f_{i}}$, i.e., the
corresponding measures on $X$. Since $\mu _{1}<<\mu _{2}$, the Radon-Nikodym
derivative $k$:~$=\frac{d\mu _{1}}{d\mu _{2}}$ is well defined. Clearly then%
\begin{equation*}
\left\Vert U\psi \right\Vert _{L^{2}(\mu _{2})}^{2}=\int_{X}\left\vert \psi
\right\vert ^{2}k\mathcal{\;}d\mu _{2}=\int_{X}\left\vert \psi \right\vert
^{2}\mathcal{\;}d\mu _{1}=\left\Vert \psi \right\Vert _{L^{2}(\mu _{1})}%
\text{.}
\end{equation*}%
As a result $W$:~$=V_{f_{2}}UV_{f_{1}}^{\ast }$ is a well defined partial
isometry in $\mathcal{H}$. For $\psi \in C(X)$, we compute 
\begin{eqnarray*}
W\rho (\psi ) &=&V_{f_{2}}UV_{f_{1}}^{\ast }\rho (\psi ) \\
&=&V_{f_{2}}UM_{\psi }V_{f_{1}}^{\ast } \\
&=&V_{f_{2}}M_{\psi }UV_{f_{1}}^{\ast } \\
&=&\rho (\psi )V_{f_{2}}UV_{f_{1}}^{\ast } \\
&=&\rho (\psi )W\text{,}
\end{eqnarray*}%
and we conclude that $W\in \mathfrak{A}^{\prime }$. The prime stands for
commutant.
\end{proof}

\begin{theorem}
\label{TheMeasTheo1}Let $N\in \mathbb{N}$, $N\geq 2$; let $\mathcal{H}$ be a
Hilbert space, and $(S_{i})$ a representation of $\mathcal{O}_{N}$ in $%
\mathcal{H}$. Let $(X,d)$ be a compact metric space, and $(\sigma
_{i})_{0\leq i<N}$ an iterated function system which is complete and
non-overlapping. Let $P$\emph{:}~$\mathcal{B}(X)\rightarrow B(\mathcal{H})$
the corresponding projection valued measure. Suppose the von Neumann algebra 
$\mathfrak{A}$ generated by $\{S_{\alpha }S_{\alpha }^{\ast }\mid k\in 
\mathbb{N}$, $\alpha \in (\mathbb{\NEG{Z}}_{N})^{k}\}$ is maximally abelian.
For $f\in \mathcal{H}$, $\left\Vert f\right\Vert =1$, set $\mu _{f}(\cdot )$%
\emph{:}~$=\left\Vert P(\cdot )f\right\Vert ^{2}$. Then the following two
conditions are equivalent\emph{:}

\noindent \emph{(}i\emph{)} $f$ is a cyclic vector for $\mathfrak{A}$.

\noindent \emph{(}ii\emph{)} $\mu _{f}\circ \sigma _{i}^{-1}<<\mu _{f}$, $%
i=0,1,\cdots ,N-1$.
\end{theorem}

\begin{proof}
We first claim that 
\begin{equation}
\mu _{f}\circ \sigma _{i}^{-1}=\mu _{S_{i}f}\text{.}  \label{TheMeasEq6}
\end{equation}

To see this, we apply (\ref{MeasEq4}) to $S_{i}f$. Then $\mu
_{S_{i}}f=\sum_{j}\mu _{S_{j}^{\ast }S_{i}f}\circ \sigma _{j}^{-1}=\mu
_{f}\circ \sigma _{i}^{-1}$ since $S_{j}^{\ast }S_{i}f=\delta _{i,j}f$. This
is the desired identity (\ref{TheMeasEq6}).

Secondly, let $i\neq j$. then naturally $S_{i}f\perp S_{j}f$. But we claim
that 
\begin{equation}
\mathcal{H}_{S_{i}f}\perp \mathcal{H}_{S_{j}f}\text{;}  \label{TheMeasEq7}
\end{equation}%
i.e., for all $A\in \mathfrak{A}$, $\left\langle S_{i}f\mid
AS_{j}f\right\rangle =0$. Since $\mathfrak{A}$ is generated by the
projections $S_{\alpha }S_{\alpha }^{\ast }$, it is enough to show that $%
S_{i}^{\ast }S_{\alpha }S_{\alpha }^{\ast }S_{j}=0$ for $\alpha =(\alpha
_{1},\cdots ,\alpha _{k})$. But $S_{i}^{\ast }S_{\alpha }S_{\alpha }^{\ast
}S_{j}=\delta _{i,\alpha _{1}}\delta _{j,\alpha _{1}}S_{\alpha _{2}}\cdots
S_{\alpha _{k}}S_{\alpha _{k}}^{\ast }\cdots S_{\alpha _{2}}^{\ast }=0$
since $i=j$. The orthogonality relation (\ref{TheMeasEq7}) follows.

We first prove (i)$\Longrightarrow $(ii); in fact we prove that $\mu
_{g}<<\mu _{f}$ for all $g\in \mathcal{H}$, if $f$ is assumed cyclic. If $f$
is cyclic, and $g\in \mathcal{H}$, $\left\Vert g\right\Vert =1$, then
clearly $\mathcal{H}_{g}\subset \mathcal{H}_{f}$. By the argument in Lemma %
\ref{TheMeasLem1}(b), we conclude that $W$:~$=V_{f}^{\ast }V_{g}$:~$%
L^{2}(\mu _{g})\rightarrow L^{2}(\mu _{f})$ commutes with the multiplication
operators. Setting $k$:$~=W(1)$, we have $\int_{X}\left\vert \psi
\right\vert ^{2}\mathcal{\;}d\mu _{g}=\int_{X}\left\vert W\psi \right\vert
^{2}\mathcal{\;}d\mu _{f}=\int \left\vert \psi \right\vert ^{2}\left\vert
k\right\vert ^{2}\mathcal{\;}d\mu _{f}$, or equivalently, $d\mu
_{g}=\left\vert k\right\vert ^{2}d\mu _{f}$. The conclusion $d\mu _{g}<<d\mu
_{f}$ follows, and $\frac{d\mu _{g}}{d\mu _{f}}=\left\vert k\right\vert ^{2}$%
.

To prove (ii)$\Longrightarrow $(i); let $i,j\in \mathbb{Z}_{N}$, and suppose 
$i\neq j.$ We saw that then $\mathcal{H}_{S_{i}f}\perp \mathcal{H}_{S_{j}f}$%
. Suppose $f$ is not cyclic. Since by (\ref{TheMeasEq6}) $\mu _{S_{i}f}=\mu
_{f}\circ \sigma _{i}^{-1}$, we get the two isometries $V_{S_{i}f}$:~$%
L^{2}(\mu _{f}\circ \sigma _{i})\rightarrow \mathcal{H}_{S_{i}f}$ with
orthogonal ranges. Let 
\begin{equation*}
k_{i}=\frac{d\mu _{f}\circ \sigma _{i}^{-1}}{d\mu _{f}}\text{, and }k_{j}=%
\frac{d\mu _{f}\circ \sigma _{j}^{-1}}{d\mu _{f}}\text{.}
\end{equation*}%
Set 
\begin{equation*}
U_{i}\psi =\psi \sqrt{k_{i}}\text{ and }U_{j}\psi =\psi \sqrt{k_{j}}\text{.}
\end{equation*}%
Then the following operator 
\begin{equation}
W\text{:~}=V_{S_{j}f}U_{j}^{\ast }U_{i}V_{S_{i}f}^{\ast }  \label{TheMeasEq8}
\end{equation}%
is well defined. It is a partial isometry in $\mathcal{H}$ with initial
space $\mathcal{H}_{S_{i}f}$ and final space $\mathcal{H}_{S_{j}f}$; i.e., $%
W^{\ast }W=\limfunc{proj}(\mathcal{H}_{S_{i}f})=p_{i}$, and $WW^{\ast }=%
\limfunc{proj}(\mathcal{H}_{S_{j}f})=p_{j}$. By the lemma, $W$ is in the
commutant of $\mathfrak{A}$. But the two projections $p_{i}$ and $p_{j}$ are
orthogonal by the lemma, i.e., $p_{i}p_{j}=0$. Relative to the decomposition 
$p_{i}\mathcal{H}\oplus p_{j}\mathcal{H}$, we now consider the following two
block matrix operators 
\begin{equation*}
\left( 
\begin{array}{cc}
0 & 0 \\ 
W & 0%
\end{array}%
\right)
\end{equation*}%
and 
\begin{equation*}
\left( 
\begin{array}{cc}
0 & W^{\ast } \\ 
0 & 0%
\end{array}%
\right) \text{;}
\end{equation*}%
and note that 
\begin{equation*}
\left( 
\begin{array}{cc}
0 & 0 \\ 
W & 0%
\end{array}%
\right) \left( 
\begin{array}{cc}
0 & W^{\ast } \\ 
0 & 0%
\end{array}%
\right) =\left( 
\begin{array}{cc}
0 & 0 \\ 
0 & p_{j}%
\end{array}%
\right) \text{,}
\end{equation*}%
while 
\begin{equation*}
\left( 
\begin{array}{cc}
0 & W^{\ast } \\ 
0 & 0%
\end{array}%
\right) \left( 
\begin{array}{cc}
0 & 0 \\ 
W & 0%
\end{array}%
\right) =\left( 
\begin{array}{cc}
p_{i} & 0 \\ 
0 & 0%
\end{array}%
\right) \text{.}
\end{equation*}%
Since the two non-commuting operators are in $\mathfrak{A}^{\prime }$, it
follows that $\mathfrak{A}^{\prime }$ is non-abelian, and as a result that $%
\mathfrak{A}$ is not maximally abelian.
\end{proof}

Two Examples: (a) Let $\mathcal{H}=L^{2}(\mathbb{T})$ where as usual $%
\mathbb{T}$ denotes the torus, equipped with Haar measure. Set $e_{n}(z)$%
\emph{:}~$=z^{n}$, $z\in \mathbb{T}$, $n\in \mathbb{Z}$, and define 
\begin{equation}
\left\{ 
\begin{array}{l}
S_{0}f(z)=f(z^{2})\;\text{for }f\in \mathcal{H}\text{, and }z\in \mathbb{T}
\\ 
S_{1}f(z)=zf(z^{2})\text{.}%
\end{array}%
\right.  \label{TheMeasEq9}
\end{equation}%
As noted in Section 1, this system is in $\limfunc{Rep}(\mathcal{O}_{2},%
\mathcal{H})$. By Theorem \ref{MeasTheo1}, there is a unique projection
valued measure $P(\cdot )$ on $\mathcal{B}([0,1))$ such that 
\begin{equation}
P\left( \left[ \frac{\alpha _{1}}{2}+\cdots +\frac{\alpha _{k}}{2^{k}}\text{%
, }\frac{\alpha _{1}}{2}+\cdots +\frac{\alpha _{k}}{2^{k}}+\frac{1}{2^{k}}%
\right) \right) =S_{\alpha }S_{\alpha }^{\ast }  \label{TheMeasEq10}
\end{equation}%
where $S_{\alpha }=S_{\alpha _{1}}\cdots S_{\alpha _{k}}$.

It is easy to check that the range of the projection $S_{\alpha }S_{\alpha
}^{\ast }$ is the closed subspace in $\mathcal{H}$ spanned by 
\begin{equation*}
\left\{ e_{n}\mid n=\alpha _{1}+2\alpha _{2}+\cdots +2^{k-1}\alpha
_{k}+2^{k}p\text{, }p\in \mathbb{Z}\right\} \text{,}
\end{equation*}%
and it follows that 
\begin{equation*}
S_{\alpha }S_{\alpha }^{\ast }e_{0}=\left\{ 
\begin{array}{l}
e_{0}\text{ if }\alpha _{1}=\alpha _{2}=\cdots =\alpha _{k}=0 \\ 
0\text{ otherwise.}%
\end{array}%
\right. 
\end{equation*}%
Hence 
\begin{equation*}
\mathcal{H}_{e_{0}}=[\mathfrak{A}_{e_{0}}]=\mathbb{C}e_{0}
\end{equation*}%
is one-dimensional, and 
\begin{equation*}
\mu _{e_{0}}(\cdot )=\left\Vert P(\cdot )e_{0}\right\Vert ^{2}=\delta _{0}
\end{equation*}%
where $\delta _{0}$ is the Dirac measure on $[0,1)$ at $x=0$. With the IFS $%
\sigma _{0}(x)=\frac{x}{2}$, $\sigma _{1}(x)=\frac{x+1}{2}$ on the
unit-interval, we get 
\begin{equation}
\left[ 
\begin{array}{c}
\mu _{e_{0}}\circ \sigma _{0}^{-1}=\delta _{0} \\ 
\mu _{e_{0}}\circ \sigma _{1}^{-1}=\delta _{\frac{1}{2}}%
\end{array}%
\right.   \label{TheMeasEq11}
\end{equation}%
making it clear that condition (ii) in Theorem \ref{TheMeasTheo1} is not
satisfied.

(b) We now modify (\ref{TheMeasEq9}) as follows:

Set 
\begin{equation}
\left\{ 
\begin{array}{l}
S_{0}f(z)=\frac{1}{\sqrt{2}}(1+z)f(z^{2})\;\text{for }f\in L^{2}(\mathbb{T})%
\text{, and }z\in \mathbb{T} \\ 
S_{1}f(z)=\frac{1}{\sqrt{2}}(1-z)f(z^{2})\text{.}%
\end{array}%
\right.  \label{TheMeasEq12}
\end{equation}%
This system $(S_{i})$ is in $\limfunc{Rep}(\mathcal{O}_{2},L^{2}(\mathbb{T}%
)) $, and $\mu _{e_{0}}(\cdot )=\left\Vert P(\cdot )e_{0}\right\Vert ^{2}=$
Lebesgue measure $dt$ on $[0,1)$, where $P(\cdot )$ is again determined by (%
\ref{TheMeasEq10}). This is the representation of$\mathcal{\ }O_{2}$ which
corresponds to the usual Haar wavelet, i.e., to 
\begin{equation}
\psi (x)=\left\{ 
\begin{array}{l}
1\text{\ if\ }0\leq x<\frac{1}{2} \\ 
-1\text{\ if\ }\frac{1}{2}\leq x<1%
\end{array}%
\right.  \label{TheMeasEq13}
\end{equation}%
and 
\begin{equation}
\psi _{j,k}(x)=2^{\frac{j}{2}}\psi (2^{j}x-k)\text{\ for }j,k\in \mathbb{Z}
\label{TheMeasEq14}
\end{equation}%
is then the standard Haar basis for $L^{2}(\mathbb{R})$; compare this with (%
\ref{TheMeasEq2})$.$ For this representation $\mathfrak{A}$ can be checked
to be maximally abelian, but it also follows from the theorem, since now the
analog of (\ref{TheMeasEq11}) is 
\begin{equation*}
\left\{ 
\begin{array}{c}
\mu _{e_{0}}\circ \sigma _{0}^{-1}=2dt\;\text{restricted to }[0,\frac{1}{2})
\\ 
\mu _{e_{0}}\circ \sigma _{1}^{-1}=2dt\;\text{restricted to }[\frac{1}{2},1)%
\text{.}%
\end{array}%
\right.
\end{equation*}%
Since $\mu _{e_{0}}=dt$ restricted to $[0,1)$, it is clear that now
condition (ii) in Theorem \ref{TheMeasTheo1} is satisfied.

\end{document}